\documentclass[11pt]{article}
\usepackage{amssymb,amsmath,amsthm,latexsym}
\title{Subspace codes in ${\rm PG(2n-1,q)}$}
\author{Antonio Cossidente\\ Dipartimento di Matematica, Informatica ed Economia\\ Universit\`a della Basilicata\\
 Contrada Macchia Romana\\ I-85100 Potenza\\ Italy\\antonio.cossidente@unibas.it\\
 Francesco Pavese\\Dipartimento di Matematica, Informatica ed Economia\\ Universit\`a della Basilicata\\
 Contrada Macchia Romana\\ I-85100 Potenza\\ Italy\\francesco.pavese@unibas.it}
\date{}

\begin{document}
\maketitle

 \newpage\noindent
 {\bf Proposed Running Head:} Subspace codes in ${\rm PG(2n-1,q)}$
 \vspace{2cm}\par\noindent {\bf Corresponding
 Author:}\\Antonio Cossidente\\
 Dipartimento di Matematica, Informatica ed Economia\\ Universit\`a della Basilicata\\
 Contrada Macchia Romana\\ I-85100 Potenza\\ Italy\\antonio.cossidente@unibas.it
\newpage
%2. Format for theorems and similar.
\newtheorem{theorem}{Theorem}[section]
\newtheorem{lemma}[theorem]{Lemma}
\newtheorem{conj}[theorem]{Conjecture}
\newtheorem{remark}[theorem]{Remark}
\newtheorem{cor}[theorem]{Corollary}
\newtheorem{prop}[theorem]{Proposition}
\newtheorem{defin}[theorem]{Definition}
\newtheorem{result}[theorem]{Result}

 \def\runningheadeven{Subspace codes in ${\rm PG(2n-1,q)$}}
\def\runningheadodd{A. Cossidente and F. Pavese}

%\makeatother
\newcommand{\Prf}{\noindent{\bf Proof}.\quad }
\renewcommand{\labelenumi}{(\arabic{enumi})}

%6. Contents of the headers

\def\cA{\mathcal A}
\def\bE{\mathbf E}
\def\bF{\mathbf F}
\def\bG{\mathbf G}
\def\bP{\mathbf P}
\def\bN{\mathbf N}
\def\bZ{\mathbf Z}
\def\bL{\mathbf L}
\def\bQ{\mathbf Q}
\def\bU{\mathbf U}
\def\bV{\mathbf V}
\def\bW{\mathbf W}
\def\bX{\mathbf X}
\def\bY{\mathbf Y}
\def\cC{\mathcal C}
\def\cD{\mathcal D}
\def\cE{\mathcal E}
\def\cF{\mathcal F}
\def\cG{\mathcal G}
\def\cH{\mathcal H}
\def\cL{\mathcal L}
\def\cM{\mathcal M}
\def\cQ{\mathcal Q}
\def\cO{\mathcal O}
\def\cP{\mathcal P}
\def\cX{\mathcal X}
\def\cY{\mathcal Y}
\def\cU{\mathcal U}
\def\cV{\mathcal V}
\def\cT{\mathcal T}
\def\cR{\mathcal R}
\def\cS{\mathcal S}
\def\cK{\mathcal K}
\def\cI{\mathcal I}
\def\PG{{\rm PG}}
\def\PGL{{\rm PGL}}
\def\GF{{\rm GF}}
\def\PSL{{\rm PSL}}
\def\GL{{\rm GL}}
\def\PGO{{\rm PGO}}

\def\ps@headings{
 \def\@oddhead{\footnotesize\rm\hfill\runningheadodd\hfill\thepage}
 \def\@evenhead{\footnotesize\rm\thepage\hfill\runningheadeven\hfill}
 \def\@oddfoot{}
 \def\@evenfoot{\@oddfoot}
}

\begin{abstract}
An $(r,M,2\delta;k)_q$ constant--dimension subspace code, $\delta >1$, is a collection $\cal C$ of $(k-1)$--dimensional projective subspaces of $\PG(r-1,q)$ such that every $(k-\delta)$--dimensional projective subspace of $\PG(r-1,q)$ is contained in at most a member of $\cal C$. Constant--dimension subspace codes gained recently lot of interest due to the work by Koetter and Kschischang \cite{KK}, where they presented an application of such codes for error-correction in random network coding.
Here a $(2n,M,4;n)_q$ constant--dimension subspace code is constructed, for every $n \ge 4$. The size of our codes is considerably larger than all known constructions so far, whenever $n > 4$. When $n=4$ a further improvement is provided by constructing an $(8,M,4;4)_q$ constant--dimension subspace code, with $M = q^{12}+q^2(q^2+1)^2(q^2+q+1)+1$.
\end{abstract}

    \section{Introduction}
Let $V$ be an $r$--dimensional vector space over $\GF(q)$, $q$ any prime power. The set $S(V)$ of all subspaces of $V$, or subspaces of the projective space $\PG(V)$, forms a metric space with respect to the {\em subspace distance} defined by $d_s(U,U')= \dim (U+U') - \dim (U \cap U')$. In the context of subspace coding theory, the main problem asks for the determination of the larger size of codes in the space $(S(V),d_s)$ ({\em subspace codes}) with given minimum distance and of course the classification of the corresponding optimal codes. Codes in the projective space and codes in the Grassmannian over a finite field referred to as subspace codes and constant--dimension codes (CDCs), respectively, have been proposed for error control in random linear network coding, see \cite{KK}.
An $(r,M,d;k)_q$ constant--dimension subspace code is a set ${\cal C}$ of $k$--subspaces of $V$, where $\vert {\cal C} \vert = M$ and minimum subspace distance $d_s({\cal C})= \min \{d_s(U,U') \; \vert \; U,U' \in {\cal C}, U \ne U' \}=d$. The maximum size of an $(r,M,d;k)_q$ constant--dimension subspace code is denoted by $\cA_q(r,d;k)$.

For general results on bounds and constructions of subspaces codes, see \cite{KSK}. More recent constructions and results can be found in \cite{ES1}, \cite{ES}, \cite{EV}, \cite{GadouleauYan}, \cite{GR}, \cite{HKK}, \cite{TR}. For a geometric approach to subspace codes see also \cite{CP}, where a connection between certain subspace codes and particular combinatorial structures is highlighted.

From a combinatorial point of view an $(r,M,2\delta;k)_q$ constant--dimension subspace code, $\delta >1$, is a collection $\cal C$ of $(k-1)$--dimensional projective subspaces of $\PG(r-1,q)$ such that every $(k-\delta)$--dimensional projective subspace of $\PG(r-1,q)$ is contained in at most a member of $\cal C$.

The set ${\cal M}_{m\times n}(q)$ of $m\times n$ matrices over the finite field $\GF(q)$ forms a metric space with respect to the {\em rank distance} defined by $d_r(A,B) = rk(A-B)$. The maximum size of a code of minimum distance $d$, $1\le d \le \min\{m,n\}$, in $({\cal M}_{m\times n}(q),d_r)$ is $q^{n(m-d+1)}$ for $m\le n$ and $q^{m(n-d+1)}$ for $m\ge n$. A code ${\cal A}\subset {\cal M}_{m\times n}(q)$ attaining this bound is said to be a $q$--ary $(m,n,k)$ {\em maximum rank distance code} ({\em MRD}), where $k=m-d+1$ for $m\le n$ and $k=n-d+1$ for $m\ge n$. A rank code $\cal A$ is called $\GF(q)$--linear if $\cal A$ is a subspace of ${\cal M}_{m\times n}(q)$. Rank metric codes were introduced by Delsarte \cite{Delsarte} and rediscovered in \cite{Gabidulin} and \cite{Roth}. Recently, these codes have found a new application in the construction of error-correcting codes for random network coding \cite{SKK}.

A {\em constant--rank code} (CRC) of constant rank $r$ in $\cM_{m \times n}(q)$ is a non--empty subset of $\cM_{m \times n}(q)$ such that all elements have rank $r$. We denote a constant--rank code with length $n$, minimum rank distance $d$, and constant--rank $r$ by $(m,n,d,r)$. The term $A(m,n,d,r)$ denotes the maximum cardinality of an $(m,n,d,r)$ constant--rank code in $\cM_{m \times n}(q)$. From \cite[Proposition 8]{GadouleauYan} we have that $A(m,n,d,r) \le \genfrac{[}{]}{0pt}{}{n}{r}_q\prod_{i=0}^{r-d}(q^m-q^i)$ and if this upper bound is attained the CRC is said to be optimal. Here $\genfrac{[}{]}{0pt}{}{n}{r}_q:= \frac{(q^n-1)\cdot \ldots \cdot(q^{n-r+1}-1)}{(q^r-1)\cdot \ldots \cdot (q-1)}$ .

In this paper we will construct a $(2n,M,4;n)_q$ constant--dimension subspace code, for every $n \ge 4$. The size of our codes is considerably larger than all known constructions so far whenever $n > 4$ (Theorem \ref{totspr}, Theorem \ref{parspr}). Our approach is completely geometric and relies on the geometry of Segre varieties. This point of view enabled us to improve (part of) the classical construction of subspaces codes arising from an MRD codes by means of certain CRCs and the geometry of a non--degenerate hyperbolic quadric of the ambient projective space.

When $n=4$, by exploring in more details the geometry of the hyperbolic quadric ${\cal Q}^+(7,q)$, a further improvement is provided by constructing an $(8,M,4;4)_q$ constant--dimension subspace code, with $M = q^{12}+q^2(q^2+1)^2(q^2+q+1)+1$. An $(8,M,4;4)_q$ constant--dimension subspace code with the same size has also been constructed in \cite{ES1} with a completely different technique. We do not know if the two constructions are equivalent but certainly both codes contain a lifted MRD code. 

In the sequel $\theta_{n,q}:= \genfrac{[}{]}{0pt}{}{n+1}{1}_q=q^n + \ldots + q + 1$ .

    \section{The geometric setting}

	\subsection{Segre varieties}
	
The {\em Segre map} may be defined as the map
$$
\sigma:\PG(n-1,q)\times\PG(n-1,q)\to \PG(n^2-1,q),
$$
taking a pair of points $x=(x_1,\dots x_n)$, $y=(y_1,\dots y_n)$ of $\PG(n-1,q)$ to their product $(x_1y_1,x_1y_2,\dots, x_ny_n)$ (the $x_iy_j$ are taken in lexicographical order). The image of the Segre map is an algebraic variety called the {\em Segre variety} and denoted by ${\cal S}_{n-1,n-1}$. The Segre variety ${\cal S}_{n-1,n-1}$ has two rulings of projective $(n-1)$--dimensional subspaces, say ${\cal R}_1$ and ${\cal R}_2$, such that two subspaces in the same ruling are disjoint, and each point of ${\cal S}_{n-1,n-1}$ is contained in exactly one member of each ruling. Also, a member of ${\cal R}_1$ meets an element of ${\cal R}_2$ in exactly one point.
From \cite[Theorem 25.5.14]{HT} certain linear sections of dimension $n(n+1)/2-1$ of ${\cal S}_{n-1,n-1}$ are Veronese varieties \cite[\S 25.1]{HT}.  For more details on Segre varieties and Veronese varieties, see \cite{HT}

	\subsection{Linear representations}	

Let $(V, k)$ be a non--degenerate formed space with associated polar space $\cP$ where $V$ is a $(d+1)$--dimensional vector space over $\GF(q^e)$ and $k$ is a sesquilinear (quadratic) form. The vector space $V$ can be considered as an $(e(d+1))$--dimensional vector space $V'$ over $\GF(q)$ via the inclusion $\GF(q) \subset \GF(q^e)$. Composition of $k$ with the trace map $T : z \in \GF(q^e) \mapsto \sum_{i=1}^{e} z^{q^i} \in \GF(q)$ provides a new form $k'$ on $V'$ and so we obtain a new formed space $(V',k')$. If our new formed space $(V',k')$ is non--degenerate, then it has an associated polar space $\cP'$. The isomorphism types and various conditions are presented in \cite{KL}, \cite{Gill}. Now each point in $\PG(d,q^e)$ corresponds to a $1$--dimensional vector space in $V$, which in turn corresponds to an $e$--dimensional vector space in $V'$, that is an $(e-1)$--dimensional projective space of $\PG(e(d + 1)-1,q)$. Extending this map from points of $\PG(d,q^e)$ to subspaces of $\PG(d,q^e)$, we obtain an injective map from subspaces of $\PG(d,q^e)$ to certain subspaces of $\PG(e(d+1)-1,q)$:
$$
 \phi: \PG(d,q^e) \rightarrow \PG(e(d+1)-1,q).
$$
The map $\phi$ is called the $\GF(q)$--{\em linear representation} of $\PG(d,q^e)$.

A {\em partial $t$--spread} of a projective space $\bP$ is a collection $\cS$ of mutually disjoint $t$--dimensional projective subspaces of $\bP$. A partial $t$--spread of $\bP$ is said to be a {\em $t$--spread} if each point of $\bP$ is contained in an element of $\bP$. The partial $t$--spread $\cS$ of $\bP$ is said to be {\em maximal}, if there is no partial $t$--spread $\cS'$ of $\bP$ containing $\cS$ as a proper subset.

The set $\cD=\{\phi(P) \; | \;\; P \in \PG(d,q^e)\}$ is an example of $(e-1)$--spread of $\PG(e(d+1)-1,q)$, called a {\it Desarguesian spread} (see \cite{Se}, Section 25). The incidence structure whose points are the elements of $\cD$ and whose lines are the $(2e-1)$--dimensional projective spaces of $\PG(e(d+1)-1,q)$ joining two distinct elements of $\cD$, is isomorphic to $\PG(d,q^e)$. One immediate consequence of the definitions is that the image of the pointset of the original polar space $\cP$ is contained in the new polar space $\cP'$ (but is not necessarily equal to it).

\subsection{A pencil of hyperbolic quadrics in $\PG(2n-1,q)$} \label{pencil}

A {\em Hermitian variety} $\cH$ of $\PG(n-1,q^2)$, is the set of absolute points for some Hermitian form defined on the underlying vector space. The variety $\cH$ is called {\em degenerate} if the corresponding Hermitian form is degenerate; else, it is called {\em non--degenerate}. Let $\cH(n-1,q^2)$ be the non--degenerate Hermitian variety of $\PG(n-1,q^2)$, $n \ge 4$ even. Then $\cH(n-1,q^2)$ has the following number of points:
$$
 \frac{(q^{n}-1)(q^{n-1}+1)}{q^2-1} .
$$
The generators of $\cH(n-1,q^2)$ are $(n-2)/2$--dimensional projective spaces and the number of generators of $\cH(n-1,q^2)$ is equal to
$$
 (q+1)(q^3+1) \cdot \ldots \cdot (q^{n-1}+1) .
$$
For further details on Hermitian varieties we refer to \cite{Se1}.

Let $\cH_1$ and $\cH_2$ be the two distinct Hermitian varieties of $\PG(n-1,q^2)$ having the following homogeneous equations
$$
 f_1 : X_1 X_{\frac{n+2}{2}}^q + \ldots + X_{\frac{n}{2}} X_{n}^q + X_1^q X_{\frac{n+2}{2}} + \ldots + X_{\frac{n}{2}}^q X_{n} = 0 ,
$$
$$
 f_2 : X_1 X_{\frac{n+2}{2}}^q + \ldots + X_{\frac{n}{2}} X_{n}^q + \omega^{q-1} ( X_1^q X_{\frac{n+2}{2}} + \ldots + X_{\frac{n}{2}}^q X_{n} ) = 0 ,
$$
respectively, where $\omega$ is a primitive element of $\GF(q^2)$. Then the Hermitian pencil $\cF$ defined by $\cH_1$ and $\cH_2$ is the set of all Hermitian varieties with equations $a f_1 + b f_2 = 0$, as $a$ and $b$ vary over the subfield $\GF(q)$, not both zero. Note that there are $q+1$ distinct Hermitian varieties in the pencil $\cF$, none of which is degenerate. The set $\cX = \cH_1 \cap \cH_2$ is called the base locus of $\cF$. Since the Hermitian varieties of a pencil cover all the points of $\PG(n-1,q^2)$, a counting argument shows that
$$
 |\cX| = \frac{(q^{n-2}+1)(q^{n}-1)}{q^2-1}
$$
and any two distinct varieties in $\cF$ intersect precisely in $\cX$.
In particular $\cX$ is a variety defined by the following equation:
$$
  X_1 X_{\frac{n+2}{2}}^q + \ldots + X_{\frac{n}{2}} X_{n}^q = 0 .
$$
Straightforward computations show that $\cX$ contains the following two $(n-2)/2$--dimensional projective spaces:
$$
 \Sigma : X_1 = \ldots = X_{\frac{n}{2}} = 0 ,  \Sigma' : X_{\frac{n+2}{2}} = \ldots = X_{n} = 0 .
$$
Also, through a point $P$ of $\Sigma$ (resp. $\Sigma'$) there pass $\theta_{\frac{n-4}{2},q^2}$ lines entirely contained in $\cX$ and these lines are contained in a generator of $\cH(n-1,q^2)$ meeting $\Sigma$ (resp. $\Sigma'$) exactly in $P$.

Let $\Pi_{r-1}$ be a $(r-1)$--dimensional projective space of $\Sigma$, $1 \le r \le (n-2)/2$, and let $\Pi_{r-1}^\perp$ be the polar space of $\Pi_{r-1}$ with respect to the unitary polarity of $\cH_1$ (or, equivalently, $\cH_2$). The intersection of $\Pi_{r-1}^\perp$ and $\Sigma'$ is a $((n-2)/2-r)$--dimensional projective space, say $\Pi'_{(n-2)/2-r}$. Note that $\langle \Pi_{r-1}, \Pi'_{(n-2)/2-r} \rangle$ is a generator of $\cH_1$ contained in $\cX$. In particular, one can see that the above construction produces
$$
\sum_{r=1}^{(n-2)/2} \genfrac{[}{]}{0pt}{}{\frac{n}{2}}{r}_{q^2}
$$
distinct generators of $\cH_1$ lying on $\cX$ and these are all the generators in common between two Hermitian varieties belonging to the pencil $\cF$ except $\Sigma$ and $\Sigma'$.

A {\em hyperbolic quadric} $\cQ^+(2n-1,q)$ of $\PG(2n-1,q)$, is the set of singular points for some non--degenerate quadratic form of hyperbolic type defined on the underlying vector space. The hyperbolic quadric $\cQ^+(2n-1,q)$ has the following number of points:
$$
 \frac{(q^{n}-1)(q^{n-1}+1)}{q-1} .
$$
The generators of $\cQ^+(2n-1,q)$ are $(n-1)$--dimensional projective spaces and the number of generators of $\cQ^+(2n-1,q)$ is equal to
$$
 2(q+1)(q^2+1) \cdot \ldots \cdot (q^{n-1}+1) .
$$
The set of all generators of the hyperbolic quadric $\cQ^+(2n-1,q)$ is divided in two distinct subsets of the same size, called {\em systems of generators} and denoted by $\cM_1$ and $\cM_2$, respectively. Let $A$ and $A'$ two distinct generators of $\cQ^+(2n-1,q)$. Then their possible intersections are projective spaces of dimension
$$
\left\{
\begin{array}{ccccccc}
0, & 2, & 4, & \ldots, & n-3 & \mbox{    if    }  & A, A' \in \cM_i, i=1,2 \\
-1, & 1, & 3, & \ldots, & n-2 & \mbox{    if    }  & A \in \cM_i, A' \in \cM_j, i,j \in \{ 1,2 \}, i \ne j
\end{array}
\right.
$$
if $n$ is odd or
$$
\left\{
\begin{array}{ccccccc}
0, & 2, & 4, & \ldots, & n-2 & \mbox{    if    }  & A \in \cM_i, A' \in \cM_j, i,j \in \{ 1,2 \}, i \ne j \\
-1, & 1, & 3, & \ldots, & n-3 & \mbox{    if    }  & A, A' \in \cM_i, i=1,2
\end{array}
\right.
$$
if $n$ is even.
For further details on hyperbolic quadrics we refer to \cite{HT}.

From \cite{KL}, if $n \ge 4$ is even, then $\phi(\cH(n-1,q^2))$ is a hyperbolic quadric $\cQ^+(2n-1,q)$ of $\PG(2n-1,q)$. In particular, points of the Hermitian variety are mapped, under the $\GF(q)$--linear representation map, to mutually disjoint lines contained in the corresponding hyperbolic quadric and covering all the points of the quadric. Now, let $\phi(\cH_i) = \cQ_i$, $i = 1,2$. Then the hyperbolic quadrics $\cQ_1$, $\cQ_2$ generate a pencil of $\PG(2n-1,q)$, say $\cF'$, containing other $q-1$ distinct hyperbolic quadrics, say $\cQ_i$, $3 \le i \le q+1$, none of which is degenerate. It turns out that the base locus of $\cF'$, say $\cX'$, consists of the
$$
 \frac{(q^{n-2}+1)(q^{n}-1)}{q-1}
$$
points covered by the lines of $\phi(\cX)$. In particular $\cX'$ contains two distinguished generators, say $S$ and $S'$, corresponding to $\Sigma$ and $\Sigma'$, respectively, that are disjoint. Hence $S$ and $S'$ belong to the same system of generators, say $\cM_1^i$ of $\cQ_i$, $1 \le i \le q+1$. Finally, if we denote by $\cG$ the set of generators meeting non--trivially both $S$ and $S'$ and belonging to each hyperbolic quadric of the pencil $\cF'$, we have that
$$
  |\cG| = \sum_{r=1}^{(n-2)/2} \genfrac{[}{]}{0pt}{}{\frac{n}{2}}{r}_{q^2} .
$$

    \section{The construction}

Let ${\cal M}_{n\times n}(q)$ be the vector space of all $n\times n$ matrices over the finite field $\GF(q)$. Let $\PG(n^2-1,q)$ be the $(n^2-1)$--dimensional projective space over $\GF(q)$ equipped with homogeneous projective coordinates $(X_1, \dots, X_{n^2})$. With the identification $a_{i+1,j} = a_{in+j}$, $0 \le i \le (n-1)$, $1 \le j \le n$, we may associate, up to a non-zero scalar factor, to a matrix $A = (a_{i,j}) \in \cM_{n \times n}(q)$ a unique point $P=(a_1, \dots, a_{n^2}) \in \PG(n^2-1,q)$, and viceversa. In this setting the Segre variety ${\cS}_{n-1,n-1}$ can be represented by all $n\times n$ matrices of rank $1$.  Let $G$ be the subgroup of $\PGL(n^2,q)$ fixing $\cS_{n-1,n-1}$, then $|G| = 2 |PGL(n,q)|$. In this context the subspace of all symmetric matrices of $\cM_{n \times n}(q)$ is represented by the $(n(n+1)/2 - 1)$--dimensional projective subspace $\Gamma$ of $\PG(n^2-1,q)$ defined by the following equations:
$$
X_{in+j} = X_{(j-1)n+i+1}, \;\;\; 0 \le i \le n-2, i+2 \le j \le n .
$$
In particular $\Gamma$ meets the Segre variety $\cS_{n-1,n-1}$ in a Veronese variety $\cV$. The subgroup of $G$ fixing $\cV$ leaves invariant a $(n(n-1)/2 - 1)$--dimensional projective subspace $\Gamma'$, which corresponds to the subspace of all skew--symmetric matrices of $\cM_{n \times n}(q)$. In particular, $\Gamma'$ is either contained in or disjoint to $\Gamma$ according as $q$ is even or odd, respectively. In any case $\Gamma'$ is disjoint from $\cS_{n-1,n-1}$.

In $\PG(n-1,q^n)$ consider a $q$-order subgeometry $\PG(n-1,q)$. Let $C \in \PGL(n,q)$ be a Singer cycle of $\PG(n-1,q)$, then $\langle C \rangle$ is a Singer cyclic group of order $\theta_{n-1,q} = (q^n-1)/(q-1)$. The group $\langle C \rangle$ partitions the points of $\PG(n-1,q^n)$ into $n$ hyperplanes and the remaining orbits are $q$-order subgeometries, see \cite{Brown}. In particular $\langle C \rangle$ fixes $n$ points in general positions and each of the $n$ fixed hyperplanes contains $n-1$ fixed points. By considering the $\GF(q)$--linear representation of $\PG(n-1,q^n)$, a point of $\PG(n-1,q^n)$ becomes a $\PG(n,q)$ that is member of a Desarguesian spread of a $\PG(n^2-1,q)$. In particular points of a $\PG(n-1,q)$ become maximal spaces of a ruling of a Segre variety $\cS_{n-1,n-1}$ of $\PG(n^2-1,q)$, see \cite{LMPT}. It follows that $\PG(n^2-1,q)$ is partitioned into $n$ $(n^2-n-1)$--dimensional projective subspaces and a certain number of Segre varieties. If $\cP$ denotes the above partition of $\PG(n^2-1,q)$, then there exists a subgroup $J$ of $G$ of order $2\theta_{n-1,q}^2$ fixing $\cP$. The group $J$ is generated by the projectivities of $\PGL(n^2,q)$ induced by $\bar{\iota}, I \otimes \bar{C}, \bar{C} \otimes I \in \GL(n^2,q)$. Here $\otimes$ denotes the Kronecker product and $C$ is induced by $\bar{C} \in \GL(n,q)$.
$$
\bar{\iota} = \left(
\begin{array}{cccc}
A_{11} & A_{21} & \dots & A_{n1}\\
A_{12} & A_{22} & \dots & A_{2n}\\
\vdots & \vdots & \ddots & \vdots\\
A_{1n} & A_{2n} & \dots & A_{nn}
\end{array}
\right) ,
$$
where $A_{ij}$ are $(n\times n)$-matrices defined as follows:
$$
A_{ij} = (a_{rs}), \;\;\; a_{rs} =
\left\{
\begin{array}{cc}
1 & (i,j) = (r,s) \\
0 & (i,j) \ne (r,s)
\end{array}
\right. .
$$
The projectivity $\iota$ induced by $\bar{\iota}$ is either an involutory homology having $\Gamma$ as axis and $\Gamma'$ as center, if $q$ is odd, or an involutory elation having $\Gamma$ as axis and $\Gamma'$ as center, if $q$ is even. Also, notice that the projectivity of $J$ induced by $\bar{C} \otimes \bar{C}$ has order $\theta_{n-1,q}$ and fixes $\cV$.

From \cite{Hup}, $\bar{C}$ is conjugate in $\GL(n,q^n)$ to the a diagonal matrix $D$
$$
D = diag (\omega, \omega^q, \dots \omega^{q^{n-1}}),
$$
for some primitive element $\omega$ of $\GF(q^n)$. In other words, there exists a matrix $E \in \GL(n,q^n)$ with $E^{-1} \bar{C} E = D$. Let $\hat{J}$ be the group generated by the projectivities of $\PGL(n^2,q^n)$ induced by $\bar{\iota}, I \otimes D, D \otimes I \in \GL(n^2,q^n)$. Since
$$
 (E \otimes E)^{-1} (I \otimes \bar{C}) (E \otimes E) = I \otimes D, \;\; (E \otimes E)^{-1} (\bar{C} \otimes I) (E \otimes E) = D \otimes I,
$$
and
$$
 (E \otimes E)^{-1} \bar{\iota} (E \otimes E) = \bar{\iota} ,
$$
it turns out that the group $\hat{J}$ fixes the $q$--order subgeometry $\Pi$ of $\PG(n^2-1,q^n)$ whose points are as follows:
$$
 (\alpha_1,\dots,\alpha_n,\alpha_n^q,\alpha_1^q,\dots,\alpha_{n-1}^q,\alpha_{n-1}^{q^2},\alpha_{n}^{q^2},\dots,\alpha_{n-2}^{q^2},\dots,\alpha_2^{q^{n-1}},\dots,\alpha_1^{q^{n-1}}) ,
$$
where $\alpha_i \in \GF(q^n)$, $1 \le i \le n$, $\prod_{i=1}^{n} \alpha_i \neq 0$, and the Segre variety $\hat{\cS}_{n-1,n-1}$ corresponding to $\cS_{n-1,n-1}$ contained in it. Furthermore the group $\hat{J}$ fixes the following $n$ $(n-1)$-dimensional projective subspaces of $\Pi$:
$$
\cX_1 = \langle U_{(a-1)n+a} \rangle, \;\; 1 \le a \le n,
$$
$$
\cX_k = \langle U_{(a_1-k)n+a_1}, U_{(n-k)n+a_2(n+1)} \rangle, \;\; k \le a_1 \le n, 1 \le a_2 \le k-1, 2 \le k \le n ,
$$
where $U_i$ denotes the point with coordinates $(0,\dots,0,1,0,\dots,0)$, with $1$ in the $i-$th position. The projectivity of $\hat{J}$ induced by $D \otimes D$ has order $\theta_{n-1,q}$ and fixes the Veronese variety $\hat{\cV} = \hat{\Gamma} \cap \hat{\cS}_{n-1,n-1}$.
In particular $\cX_1$ is contained in $\hat{\Gamma}$ and the involution $\iota$ fixes $\cX_1$ pointwise and interchanges $\cX_k$ with $\cX_{n-k+2}$, $2 \le k \le n$. Then the involution $\iota$ fixes the $(n^2-n-1)$--dimensional projective subspace $\hat{\cY} = \langle \cX_k \rangle$, $2 \le k \le n$. It follows that the center of $\iota$, $\hat{\Gamma}'$, must be contained in $\hat{\cY}$. We have proved the following result.

\begin{prop}
There exists an $(n^2-n-1)$--dimensional projective space $\cY$ that is disjoint from $\cS_{n-1,n-1}$ and contains $\Gamma'$.
\end{prop}

We denote by $\cA$ the set consisting of $q^{n(n-1)/2}$ matrices corresponding to the points of $\Gamma'$ (together with the zero matrix). Since $\cY$ is disjoint from the Segre variety $\cS_{n-1,n-1}$, the set $\cM$, consisting of the $q^{n^2-n}$ matrices corresponding to the points of $\cY$ (together with the zero matrix), form a linear $(n,n,n-1)$ MRD code.

Let $A$ be a $n \times n$ matrix over $\GF(q)$, and let $I_n$ be the $n \times n$ identity matrix. The rows of the $n \times 2n$ matrix $(I_n | A)$ can be viewed as points in general position of an $(n-1)$-dimensional projective space of $\PG(2n-1,q)$. This subspace is denoted by $L(A)$. From \cite{SKK}, a $q$-ary $(n,n,n-1)$ MRD lifts to a $q$-ary $(2n,q^{n^2-n},4;n)$ constant--dimension subspace code. A constant--dimension code such that all its codewords are lifted codewords of an MRD code is called a {\em lifted MRD code}. Let $\cL_1 = \{ L(A) | A \in \cM \}$ be the constant--dimension code obtained by lifting the $(n,n,n-1)$ MRD code contructed above. Then $\cL_1$ consists of $(n-1)$--dimensional projective spaces mutually intersecting in at most an $(n-3)$--dimensional projective space. In particular, members of $\cL_1$ are disjoint from the special $(n-1)$-dimensional projective space $S=\langle U_{n+1}, \dots, U_{2n} \rangle$ and therefore every $(n-2)$--dimensional projective space covered by an element of $\cL_1$ is disjoint from $S$. Moreover, from \cite[Lemma 6]{HKK}, every $(n-2)$-dimensional projective space in $\PG(2n-1,q)$ disjoint from $S$ is covered by a member of $\cL_1$ exactly once.

From \cite{Gabidulin} it is known that a linear $(n,n,n-1)$ MRD code contains an $(n,n,2,r)$ CRC of size
$$
 \genfrac{[}{]}{0pt}{}{n}{r}_q \sum_{j=2}^{r} (-1)^{(r-j)} \genfrac{[}{]}{0pt}{}{r}{j}_q q^{\genfrac{(}{)}{0pt}{}{r-j}{2}}(q^{n(j-1)}-1) .
$$
Let $\cC_r$ denotes the $(n,n,2,r)$ CRC contained in $\cY$. Let $A$ be an element of $\cC_r$, $2 \le r \le (n-2)$. Again, the rows of the $n \times 2n$ matrix $(A | I_n)$ can be viewed as points in general position of an $(n-1)$-dimensional projective space of $\PG(2n-1,q)$. This subspace is denoted by $L'(A)$. The subspace $L'(A)$ is disjoint from the special $(n-1)$-dimensional projective space $S'=\langle U_{1}, \dots, U_{n} \rangle$ and meets $S$ in a $(n-r-1)$--dimensional projective space. It follows that every $(n-2)$--dimensional projective space contained in $L'(A)$ meets $S$ in at least a point and is disjoint from $S'$. Let $\cL_r = \{ L'(A) | A \in \cC_r \}$ be the constant--dimension code obtained by lifting the $(n,n,2,r)$ CRC codes $\cC_r$, $2 \le r \le (n-2)$ constructed above. If $A_1 \in \cC_{r_1}$, $A_2 \in \cC_{r_2}$, then $L'(A_1)$ meets $L'(A_2)$ in at most in $(n-3)$--dimensional projective space. Then we have the following result:
\begin{prop}
The set $\bigcup_{i=1}^{n-2} \cL_i$ is a $(2n,M,4;n)_q$ constant--dimension subspace code, where
$$
M = q^{n^2-n} + \sum_{r=2}^{n-2} \genfrac{[}{]}{0pt}{}{n}{r}_q \sum_{j=2}^{r} (-1)^{(r-j)} \genfrac{[}{]}{0pt}{}{r}{j}_q q^{\genfrac{(}{)}{0pt}{}{r-j}{2}}(q^{n(j-1)}-1) .
$$
\end{prop}
	
Now, we introduce the non--degenerate hyperbolic quadric $\cQ$ of $\PG(2n-1,q)$ having the following equation:
$$
    X_1 X_{2n} + X_2 X_{2n-1} + \ldots + X_n X_{n+1} = 0 .
$$  	
The subspaces $S$ and $S'$ are maximals of $\cQ$ that belong to the same or to different systems, according as $n$ is even or odd, respectively. Let $\cM_1$ be the system of maximals of $\cQ$ containing $S$ and let $D(X)$ and $I(X)$ denote the set of maximals in $\cM_1$ disjoint from $X$ or meeting non--trivially $X$, respectively. Let $A$ be a skew--symmetric matrix in $\cM_{n \times n}(q)$, then it is not difficult to see that $L(A)$ (resp. $L'(A)$) is a maximal of $\cQ$ disjoint from $S$ (resp. $S'$). Since the number of maximals of $\cQ$ disjoint from $S$ equals $q^{n(n-1)/2}$ \cite[p. 175 Ex. 11.4]{Taylor}, we have that each such a maximal is of the form $L(A)$, for some $A \in \cA$.

\subsection{$n$ even}

Assume that $n$ is even. In this case we have that
$$
\cM_1 = D(S) \cup (D(S') \cap I(S)) \cup (I(S) \cap I(S'))
$$
and
$$
 |D(S)| = q^{\frac{n(n-1)}{2}} .
$$
On the other hand, a maximal $L'(A)$ in $D(S')$ is disjoint from $S$ if and only if $A$ is a skew--symmetric matrix of rank $n$. From \cite{Lewis}, the number of skew--symmetric matrices of rank $n$ is equal to
$$
  q^{\frac{n(n-2)}{4}}(q^{n-1}-1) (q^{n-3}-1) \cdot \ldots \cdot (q-1) = q^{\frac{n(n-2)}{4}} \prod_{i=0}^{\frac{n-2}{2}} (q^{2i+1}-1) .
$$
Therefore, we have that
$$
 |D(S') \cap I(S)| = q^{\frac{n(n-1)}{2}} - q^{\frac{n(n-2)}{4}} \prod_{i=0}^{\frac{n-2}{2}} (q^{2i+1}-1)
$$
and
$$
 |I(S) \cap I(S')| = |\cM_1| -  2 q^{\frac{n(n-1)}{2}} + q^{\frac{n(n-2)}{4}} \prod_{i=0}^{\frac{n-2}{2}} (q^{2i+1}-1) .
$$
Notice that both $D(S)$ and $D(S') \cap I(S)$ are contained in $\bigcup_{i=1}^{n-2} \cL_i$, whereas $I(S) \cap I(S')$ is disjoint from $\bigcup_{i=1}^{n-2} \cL_i$. Then it turns out that $( \bigcup_{i=1}^{n-2} \cL_i ) \cup (I(S) \cap I(S'))$  is a set of $(n-1)$--dimensional projective spaces mutually intersecting in at most an $(n-3)$--dimensional projective space of size
$$
 q^{n^2-n} + \sum_{r=2}^{n-2} \genfrac{[}{]}{0pt}{}{n}{r}_q \sum_{j=2}^{r} (-1)^{(r-j)} \genfrac{[}{]}{0pt}{}{r}{j}_q q^{\genfrac{(}{)}{0pt}{}{r-j}{2}}(q^{n(j-1)}-1) +
$$
$$
+ \prod_{i=1}^{n-1} (q^i+1) -  2 q^{\frac{n(n-1)}{2}} + q^{\frac{n(n-2)}{4}} \prod_{i=0}^{\frac{n-2}{2}} (q^{2i+1}-1) .
$$

In this case, from Section (\ref{pencil}), there exists a pencil $\cF'$ comprising $q$ hyperbolic quadrics $\cQ_i$, $2 \le i \le q+1$ of $\PG(2n-1,q)$ distinct from $\cQ$. Let $I_i(X)$ denote the set of maximals in $\cM^i_1$ meeting non--trivially $X$, $2 \le i \le (q+1)$ and let $\cG = \bigcap_{i=2}^{q+1}(I_i(S) \cap I_i(S')) \cap (I(S) \cap I(S'))$. Then, from Section (\ref{pencil}), we have that
$$
 |\bigcup_{i = 2}^{q+1} (I_i(S) \cap I_i(S'))| = q (|I(S) \cap I(S')| - |\cG|) =
$$
$$
 q \left(|\cM_1| -  2 q^{\frac{n(n-1)}{2}} + q^{\frac{n(n-2)}{4}} \prod_{i=0}^{\frac{n-2}{2}} (q^{2i+1}-1) - \sum_{r=1}^{\frac{n-2}{2}} \genfrac{[}{]}{0pt}{}{\frac{n}{2}}{r}_{q^2} \right) .
$$
It follows that $( \bigcup_{i=1}^{n-2} \cL_i ) \cup ( \bigcup_{i = 2}^{q+1} (I_i(S) \cap I_i(S')) ) \cup (I(S) \cap I(S'))$ is a set of $(n-1)$--dimensional projective spaces mutually intersecting in at most an $(n-3)$--dimensional projective space of size
$$
 q^{n^2-n} - 2 (q+1) q^{\frac{n(n-1)}{2}} + \sum_{r=2}^{n-2} \genfrac{[}{]}{0pt}{}{n}{r}_q \sum_{j=2}^{r} (-1)^{(r-j)} \genfrac{[}{]}{0pt}{}{r}{j}_q q^{\genfrac{(}{)}{0pt}{}{r-j}{2}}(q^{n(j-1)}-1) +
$$
$$
 + (q+1) \left( \prod_{i=1}^{n-1} (q^i+1) + q^{\frac{n(n-2)}{4}} \prod_{i=0}^{\frac{n-2}{2}} (q^{2i+1}-1) \right) - q \sum_{r=1}^{\frac{n-2}{2}} \genfrac{[}{]}{0pt}{}{\frac{n}{2}}{r}_{q^2} .
$$

The set $\cG$ contains a subset $\cD$ consisting of $\theta_{(n-2)/2,q^2}$ generators belonging to each hyperbolic quadric of the pencil $\cF'$ such that every element in $\cD$ meets $S$ in a line and $S'$ in an $(n-3)$--dimensional projective space and the set $\cD_S = \{ A \cap S \; | \;\; A \in \cD \}$ is a Desarguesian line--spread of $S$. In other words $\cD_{S} = \{ \phi(P) \; | \;\; P \in \Sigma \}$. On the other hand, the set $\cD_{S'} = \{ A \cap S' \; | \;\; A \in \cD \}$ is a set of $(n-3)$--dimensional projective space mutually intersecting in an $(n-5)$--dimensional projective space. In particular for a fixed line $\ell \in \cD_S$ there exists a unique element in $\cD_{S'}$, say $A_{\ell}$, such that $\langle \ell, A_{\ell} \rangle$ is in $\cD$, and viceversa. Furthermore, if $\ell \in \cD_{S}$ and $B \in \cD_{S'} \setminus \{ A_{\ell} \}$, then $\langle \ell, B \rangle$ is an $(n-1)$--dimensional projective space meeting a hyperbolic quadric of the pencil $\cF'$ in a cone having as vertex $A_{\ell} \cap B$ and as base a $\cQ^+(3,q)$ containing $\ell$. Notice that such a cone meets a generator of a hyperbolic quadric of the pencil $\cF'$ in at most an $(n-3)$--dimensional projective space. Let $\cD'$ be the set of $(n-1)$--dimensional projective spaces of the form $\langle \ell, B \rangle$, where $\ell \in \cD_{S}$ and $B \in \cD_{S'} \setminus \{ A_{\ell} \}$. Then $\cD'$ is disjoint from $\cD$. Also $|\cD'| = \theta_{(n-2)/2,q^2}(\theta_{(n-2)/2,q^2}-1)$.
From the discussion above, we have that $( \bigcup_{i=1}^{n-2} \cL_i ) \cup ( \bigcup_{i = 2}^{q+1} (I_i(S) \cap I_i(S')) ) \cup (I(S) \cap I(S')) \cup \cD' \cup \{ S \}$ is a set of $(n-1)$--dimensional projective spaces mutually intersecting in at most an $(n-3)$--dimensional projective space. We have proved the following result.

\begin{theorem}\label{totspr}
If $n$ is even, there exists a $(2n,M,4;n)_q$ constant--dimension subspace code, where
$$
M = q^{n^2-n} - 2 (q+1) q^{\frac{n(n-1)}{2}} + \sum_{r=2}^{n-2} \genfrac{[}{]}{0pt}{}{n}{r}_q \sum_{j=2}^{r} (-1)^{(r-j)} \genfrac{[}{]}{0pt}{}{r}{j}_q q^{\genfrac{(}{)}{0pt}{}{r-j}{2}}(q^{n(j-1)}-1) +
$$
$$
 + (q+1) \left( \prod_{i=1}^{n-1} (q^i+1) + q^{\frac{n(n-2)}{4}} \prod_{i=0}^{\frac{n-2}{2}} (q^{2i+1}-1) \right) - q \sum_{r=1}^{\frac{n-2}{2}} \genfrac{[}{]}{0pt}{}{\frac{n}{2}}{r}_{q^2} + \genfrac{[}{]}{0pt}{}{\frac{n}{2}}{1}_{q^2} \left( \genfrac{[}{]}{0pt}{}{\frac{n}{2}}{1}_{q^2} - 1 \right) + 1 .
$$
\end{theorem}

\subsection{$n$ odd}

Assume that $n$ is odd. In this case
$$
\cM_1 = (D(S') \cap I(S)) \cup (I(S) \cap I(S'))
$$
and
$$
 |D(S)| = 0 , \;\;\; |D(S') \cap I(S)| = |D(S')| = q^{\frac{n(n-1)}{2}} .
$$
On the other hand, a maximal $L'(A)$ in $D(S')$ is not in $\bigcup_{i=1}^{n-2} \cL_i$ if and only if $A$ is a skew--symmetric matrix of rank $n-1$, i.e.,  $L'(A)$ meets $S$ in a point. From \cite{Lewis}, the number of skew--symmetric matrices of rank $n-1$ is equal to
$$
  q^{\frac{(n-1)(n-3)}{4}}(q^{n}-1) (q^{n-2}-1) \cdot \ldots \cdot (q^3-1) = q^{\frac{(n-1)(n-3)}{4}} \prod_{i=1}^{\frac{n-1}{2}} (q^{2i+1}-1) .
$$
Therefore, if we denote by $\cI$ the subset of $D(S')$ consisting of maximal intersecting $S$ in exactly a point, we have that
$$
 |\cI| = q^{\frac{(n-1)(n-3)}{4}} \prod_{i=1}^{\frac{n-1}{2}} (q^{2i+1}-1)
$$
and
$$
 |I(S) \cap I(S')| = |\cM_1| -  q^{\frac{n(n-1)}{2}} .
$$
Notice that $\{ L(A) | A \in \cA\} \subseteq \cL_1$. Then, if $\cL'_1 = \cL_1 \setminus \{ L(A) | A \in \cA\}$, then $\cL'_1 \cup ( \bigcup_{i=2}^{n-2} \cL_i ) \cup \cI \cup (I(S) \cap I(S')) \cup \{ S \}$ is a set of $(n-1)$--dimensional projective spaces mutually intersecting in at most an $(n-3)$--dimensional projective space of size
$$
 q^{n^2-n} + \sum_{r=2}^{n-2} \genfrac{[}{]}{0pt}{}{n}{r}\sum_{j=2}^{r} (-1)^{(r-j)} \genfrac{[}{]}{0pt}{}{r}{j}q^{\genfrac{(}{)}{0pt}{}{r-j}{2}}(q^{n(j-1)}-1) +
$$
$$
 + \prod_{i=1}^{n-1} (q^i+1) - 2 q^{\frac{n(n-1)}{2}} + q^{\frac{(n-1)(n-3)}{4}} \prod_{i=1}^{\frac{n-1}{2}} (q^{2i+1}-1) + 1 .
$$

From \cite[Theorem 4.6]{B} a partial $1$--spread of $\PG(n-1,q)$, $n \ge 5$ odd, has size $y:=q^{n-2}+q^{n-4}+\dots+q^3+1$ and actually examples of this size exist. Arguing as in the last part of the previous paragraph we prove the following result.

\begin{theorem}\label{parspr}
If $n$ is odd, there exists a $(2n,M,4;n)_q$ constant--dimension subspace code, where
$$
M=q^{n^2-n} + \sum_{r=2}^{n-2} \genfrac{[}{]}{0pt}{}{n}{r}\sum_{j=2}^{r} (-1)^{(r-j)} \genfrac{[}{]}{0pt}{}{r}{j}q^{\genfrac{(}{)}{0pt}{}{r-j}{2}}(q^{n(j-1)}-1) +
$$
$$
 + \prod_{i=1}^{n-1} (q^i+1) - 2 q^{\frac{n(n-1)}{2}} + q^{\frac{(n-1)(n-3)}{4}} \prod_{i=1}^{\frac{n-1}{2}} (q^{2i+1}-1)+y(y-1) + 1 .
$$
\end{theorem}

\section{The case of $\PG(7,q)$}

In this section we will improve, in the case $n=4$, the result established in Theorem \ref{totspr} by considering some more suitable projective $3$--spaces (solids).

In this case $S$ and $S'$ are generators of ${\cal Q}^+(7,q)$ belonging to the same system.
Here, $\cD$ consists of $q^2+1$ generators belonging to each hyperbolic quadric of the pencil $\cF'$ such that every element in $\cD$ meets $S$ and $S'$ in a projective line.  It follows that $\cD_{S} = \{ A \cap S \; | \;\; A \in \cD \}$ and $\cD_{S'} = \{ A \cap S' \; | \;\; A \in \cD \}$ are both Desarguesian line--spreads of $S$ and $S'$, respectively. In other words $\cD_{S} = \{ \phi(P) \; | \;\; P \in \Sigma \}$ and $\cD_{S'} = \{ \phi(P) \; | \;\; P \in \Sigma' \}$.
Let $r'$ be a line of $S'$. Then, $r'^\perp$ (here $\perp$ denotes the orthogonal polarity of $\PG(7,q)$ induced by ${\cal Q}^+(7,q)$) meets $S$ in a line $r$. If $r'$ belongs to $\cD_{S'}$, then $r$ belongs to $\cD_{S}$. Assume that $r'$ does not belong to $\cD_{S'}$. Of course, $r'$ meets $q+1$ lines $l_1',\dots,l_{q+1}'$ of $\cD_{S'}$ and $r$ meets $q+1$ lines $l_1,\dots,l_{q+1}$ of $\cD_{S}$.
The subgroup of the orthogonal group $\PGO^+(8,q)$ fixing $\cQ^+(7,q)$ and stabilizing both $S$ and $S'$ (but that does not interchange them) is isomorphic to $\PGL(4,q)$ (which in turn is isomorphic to a subgroup of index two of $\PGO^+(6,q)$).
Under the Klein correspondence between lines of $S$ and points of the Klein quadric $\cal K$, the lines of $\cD_{S}$ are mapped to a $3$--dimensional elliptic quadric $\cal E$ embedded in $\cal K$ and the lines $l_1,\dots,l_{q+1}$ are mapped to a conic section $\cal C$ of $\cal E$, see \cite{Hir}. Also, there exists a subgroup $H'$ of the orthogonal group $\PGO^+(6,q)$ fixing $\cK$, isomorphic to $C_{q+1} \times \PGL(2,q^2)$, stabilizing $\cE$ and permuting in a single orbit the remaining points of $\cK$. It follows that there exists a group $H$ in the orthogonal group $\PGO^+(8,q)$ corresponding to $H'$,  stabilizing $\cQ^+(7,q)$ and fixing both $S$, $S'$, their line--spreads $\cD(S)$, $\cD(S')$ and permuting in a single orbit the remaining lines of $S$ (respectively $S'$). In this setting the line $r$ corresponds, under the Klein correspondence, to a point $P \in \cC^{\perp_{\cK}}$ (here $\perp_{\cK}$ denotes the orthogonal polarity of $\PG(5,q)$ induced by $\cK$). Let $H'_P$ be the stabilizer of $P$ in $H'$. Then $|H'_P| = |\PGL(2,q)|$. On the other hand, $H'_{\cal C}$, the stabilizer of $\cal C$ in $H'$, is contained in $H'_P$ and contains a subgroup isomorphic to $\PGL(2,q)$. It follows that $H'_P = H'_{\cC} \simeq \PGL(2,q)$. The group $H'_{\cal C}$ has $q(q-1)/2$ orbits of size $q^2-q$. Each of them together with $\cal C$ gives rise to an elliptic quadric of $\cal K$ on $\cal C$ and these are all the elliptic quadrics of $\cal K$ on $\cal C$. Let ${\cal E}'$ be one of the above orbits of $H'_{\cal C}$ of size $q^2-q$ disjoint from $\cal E$. Let $L_{{\cal E}'}$ be the set of lines of $S$ corresponding to ${\cal E}'$. Let $Y$ denotes the solid generated by $r'$ and a line of $L_{{\cal E}'}$ and consider the orbit $Y^H$ of $Y$ under the action of the group $H$. Since the lines in $L_{\cE'}$ are mutually disjoint, then two distinct solids in $Y^H$ containing $r'$ have in common exactly the line $r'$. Let $l$ be a line of $L_{\cE'}$. Under the Klein correspondence, the line $l$ corresponds to a point $P' \in \cE'$. Notice that $P'^{\perp_{\cK}}$ meets $\cE$ in a conic, say $\cC'$, that is necessarily disjoint from $\cC$. Assume on the contrary that there exists a point in common between $\cC$ and $\cC'$, say $Q$. Then the line $P'Q$ is entirely contained in $\cK$. Also, $P'Q \subset \cE' = \langle P', \cC \rangle \cap \cK$, contradicting the fact that $\cE'$ is a $3$--dimensional elliptic quadric (and so does not contain lines). Now, we claim that the solid $\langle P, \cC' \rangle$ meets $\cK$ in a $3$--dimensional elliptic quadric. Indeed, otherwise, there would be a line entirely contained in $\cK$ and passing through $P$. But such a line would contain a point of $\cC'$, that clearly is a contradiction, since $P \in \cC^{\perp_{\cK}}$ and $\cC'$ is disjoint from $\cC$. It follows that if $H_l$ denotes the stabilizer of $l$ in $H$, then $r^{H_l}$ contain $q^2-q$ mutually disjoint lines. Therefore $r'^{H_l}$ contain $q^2-q$ mutually disjoint lines and two solids in $Y^H$ containing $l$ have in common exactly the line $l$. Then $Y^H$ is a set of solids mutually intersecting in at most a line. The set $Y^H$ contains $(q^2-q)(q^2+q)(q^2+1) = q^6-q^2$ solids. Notice that none of the solids in $Y^H$ is a generator of $\cQ^+(7,q)$ or of a quadric of the pencil $\cF'$. Finally, assume that a solid $T$ in $Y^H$ generated by a line $l\in S$ and a line $r\in S'$ contains a plane $\pi$ that is entirely contained in $\cQ^+(7,q)$ or in a quadric of the pencil $\cF'$. Then, $\pi$ would meet $l'$ in a point $U$ and hence $T$ would meet $S'$ in a line through $U$ that is not the case.  We have proved the following result.

\begin{theorem} \label{triality}
There exists an $(8,M,4;4)_q$ constant--dimension subspace code, where
$$
M = q^{12}+q^2(q^2+1)^2(q^2+q+1)+1 .
$$
\end{theorem}

\begin{cor}
$$
 \cA_q(8,4;4) \ge q^{12}+q^2(q^2+1)^2(q^2+q+1)+1 .
$$
\end{cor}

\begin{remark}
{\rm The result of Theorem \ref{triality} was obtained with different techniques in \cite{ES1}, where the authors, among other interesting results, proved that $q^{12}+q^2(q^2+1)^2(q^2+q+1)+1$ is also the maximum size of an $(8,M,4;4)_q$ constant--dimension subspace code containing a lifted MRD code.}
\end{remark}


\begin{thebibliography}{10}

\bibitem{B} A. Beutelspacher, Partial spreads in finite projective spaces and partial designs, {\em Math. Z.} 145 (1975), no. 3, 211-229.

\bibitem{Brown} J. M. N. Brown, Partitioning the complement of a simplex in $\PG(e,q^{d+1})$ into copies of $\PG(d,q)$, {\em J. Geom.} 33 (1988), no. 1-2. 11-16.

\bibitem{CP} A. Cossidente, F. Pavese, On subspace codes, {\em Des. Codes Cryptogr.} DOI 10.1007/s10623-014-0018-6.

\bibitem{Delsarte} P. Delsarte, Bilinear forms over a finite field, with applications to coding theory, {\em  J. Combin. Theory Ser. A }, 25 (1978)  226-241.

\bibitem{ES1} T. Etzion, N. Silberstein, Codes and Designs Related to Lifted MRD Codes, {\em IEEE Trans. Inform. Theory} 59 (2013), no. 2, 1004-1017.

\bibitem{ES} T. Etzion, N. Silberstein, Error-correcting codes in projective spaces via rank-metric codes and Ferrers diagrams, {\em IEEE Trans. Inform. Theory} 55 (2009), no. 7, 2909-2919.

\bibitem{EV} T. Etzion, A. Vardy, Error-correcting codes in projective space, {\em  IEEE Trans. Inform. Theory} 57 (2011), no. 2, 1165-1173.

\bibitem{Gabidulin} E. M. Gabidulin, Theory of codes with maximum rank distance, {\em Problems of Information Transmission} 21 (1985), 1-12.

\bibitem{GadouleauYan} M. Gadouleau, Z. Yan, Constant-rank codes and their connection to constant-dimension codes, {\em IEEE Trans. Inform. Theory} 56 (2010), no. 7, 3207-3216.

\bibitem{Gill} N. Gill, Polar spaces and embeddings of classical groups, {\em New Zealand J. Math.} 36 (2007), 175-184.

\bibitem{GR} E. Gorla, A. Ravagnani, Subspace codes from Ferrers diagrams, preprint (arXiiv:1405.2736).

\bibitem{Hir}  J. W. P. Hirschfeld, {\em Finite projective spaces of three dimensions}, Oxford Mathematical Monographs, Oxford Science Publications, The Clarendon Press, Oxford University Press, New York, 1985.

\bibitem{HT} J. W. P. Hirschfeld, J. A. Thas, {\em General Galois Geometries}, Oxford Mathematical Monographs, Oxford Science Publications, The Clarendon Press, Oxford University Press, New York, 1991.

\bibitem{HKK} T. Honold, M. Kiermaier, S. Kurz, Optimal binary subspace codes of length $6$, constant dimension $3$ and minimum distance $4$, preprint (arXiiv:1311.0464v1).

\bibitem{Hup} B. Huppert, {\em Endliche Gruppen I}, Die Grundlehren der Mathematischen Wissenschaften, Springer-Verlag, Berlin-New York, 1967.

\bibitem{KSK} A. Khaleghi, D. Silva, F. R. Kschischang, Subspace codes, {\em Cryptography and coding}, 1-21, Lecture Notes in Comput. Sci., Springer, Berlin, 2009.

\bibitem{KL} P. Kleidman, M. Liebeck, {\em The subgroup structure of the finite classical groups}, London Mathematical Society Lecture Note Series, vol. 129 Cambridge University Press, Cambridge, 1990.

\bibitem{KK} R. Koetter, F. R. Kschischang, Coding for errors and erasures in random network coding, {\em IEEE Trans. Inform. Theory}, 54 (8), 3579-3591, Aug. 2008.

\bibitem{Lewis} J. B. Lewis, R. I. Liu, A. H. Morales, G. Panova, S. V. Sam, Y. X. Zhang, Matrices with restricted entries and $q$--analogues of permutations, {\em J. Comb.} 2 (2011), no. 3, 355-395.

\bibitem{LMPT} G. Lunardon, G. Marino, O. Polverino, R. Trombetti, Maximum scattered linear sets of pseudoregulus type and the Segre Variety ${\cal S}_{n,n}$, {\em J. Algebr. Comb.}, 39 (2014), 807-831.

\bibitem{Roth} R. M. Roth, Maximum-rank array codes and their application to crisscross error correction, {\em IEEE Trans. Inform. Theory}, 37 (1991), 328-336.

\bibitem{Se1} B. Segre, Forme e geometrie hermitiane, con particolare riguardo al caso finito,  {\em Ann. Mat. Pura Appl.} 70 (4) (1965),  1-201.

\bibitem{Se} B. Segre, Teoria di Galois, fibrazioni proiettive e geometrie non desarguesiane, {\em Ann. Mat. Pura Appl.} 64 (1964), 1-76.

\bibitem{SKK} D. Silva, F. R. Kschischang, R. Koetter, A rank-metric approach to error control in random network coding, {\em IEEE Trans. Inform. Theory}, vol. 54, pp. 3951-3967, September 2008.

\bibitem{TR} A.-L. Trautmann, J. Rosenthal, New improvements on the echelon-Ferrers construction, in proc. of {\em Int. Symp. on Math. Theory of Networks and Systems}, 405-408, July 2010.

\bibitem{Taylor} D. E. Taylor, {\em The geometry of the classical groups}, Sigma Series in Pure Mathematics, 9, Heldermann Verlag, Berlin, 1992.

    \end{thebibliography}
    \end{document}